\documentclass[11pt,twoside]{amsart}

\usepackage[english]{babel}

\usepackage{amsmath,amsfonts,amsthm}
\usepackage{graphicx}
\usepackage{epsfig}

\def\ds{\displaystyle}
\def\eps{\varepsilon}
\def\R{{\mathbb R}}

\def\N{{\mathbb N}}
\def\e{{\bf e}}
\def\o{\overline}

\newtheorem{theo}{Theorem}
\newtheorem{prop}{Proposition}
\newtheorem{cor}{Corollary}
\newtheorem{lem}{Lemma}

\author[N. Enriquez]{Nathana\"el ENRIQUEZ}
\address{Laboratoire de Probabilit\'es de Paris 6, 4
place Jussieu, 75252 Paris cedex 05}
\email{enriquez@ccr.jussieu.fr}

\title[Az\'ema martingales]{An invariance principle for Az\'ema martingales.}

\begin{document}
\maketitle
{\bf Abstract}:
{\small An invariance principle for Az\'ema martingales is presented as well as 
 a new device to construct solutions of Emery's structure equations.}

\section{Introduction}

 In \cite{M1}, Meyer raised the problem of finding ``normal" martingales, i.e.
martingales $M_t$ on a filtered probability space $(\Omega, {\mathcal
F}, {\mathcal F}_t, P)$ satisfying
$\langle M,M\rangle_t=t$, which would enjoy the chaotic representation property.
This property, which is stronger than the previsible representation property,
requires the direct sum of the chaoses of $M$ to be equal to the whole space
$L^2(\Omega)$. The only known examples of such martingales were the Brownian
motion and the compensated Poisson process until Emery found in \cite {E1} a
way to construct a whole family of examples, by introducing the beautiful device
of structure equations. Among these examples, some of them, which are
self-similar, appeared to be the generalization of two martingales considered
earlier, for different purposes. One had been  considered by Az\'ema in
\cite{A1},  and the other by Parthasarathy in an unpublished paper \cite{P1}.
The whole class of these self-similar martingales is now usually called Az\'ema 
martingales. However, their behavior remains a little mysterious,
especially, in the neigbourhood of 0. We refer the reader to the two lecture
notes of Yor \cite{Y}, and of Mansuy and Yor \cite{MY}, where a chapter is
devoted to the subject. 

Inspired by some elementary remarks on
renewal processes, we present an invariance principle for Az\'ema martingales,
and propose a new device to construct solutions of Emery's structure equations.
The principle of these approximations turns out to be quite unusual and differs
from the truncation method one can find in dealing with L\'evy processes, or from
the discretization scheme already used in the subject by Meyer in \cite{M2} to
construct solutions of Emery's structure equations. It is mainly based on the
introduction of some randomness in the size of the jumps together with keeping
the normal martingale property.

\section{The basic example: the first Az\'ema martingale}

 The aim of this section is to introduce, in a natural way, the device
we use in the next section to approximate general Az\'ema martingales.  

Let us start with the first Az\'ema martingale. It was obtained by Az\'ema in
\cite{A1}, by projecting  a Brownian motion $B_t$ starting from $0$, on the
filtration of
$sign\,(B)$. This projection yields the following normal martingale 
$$X_t={\sqrt{\pi\over2}} sign\, (B_t)\sqrt {(t-g_t)} $$
where $g_t=\sup\{ s : s\leq t,\, B_s=0\}$. We refer the reader to Chapter IV of
Protter's book \cite{Pr} for a comprehensive discussion about the martingale
property of $X_t$. A clear justification of the projection property is
given in an article of Az\'ema and Yor
\cite{AY}: they use an elegant path decomposition of the Brownian motion
on $[0,t]$ involving the Brownian meander, seen as the renormalisation of the
last incomplete excursion.

Now, it is well known, from P. L\'evy, that
$2X_t^2/\pi t$ follows the arcsine law. On the other hand, there is an old
result of Dynkin \cite{D} presented also in Feller's book (\cite{Feller} Chap.
XIV.3), obtaining arcsine laws from asymptotic waiting times of renewal processes
having non-integrable inter-arrival times.  Now, the question is:  can we
associate a martingale to any given renewal process, exactly like the Az\'ema
martingale is associated to the zeroes of the Brownian motion?  

The following elementary proposition gives a positive answer to this question:

\begin{prop}
Let $X_1, X_2,...$ be an iid sequence of positive random variables admitting a
density function. We denote by
$S_n:=\sum_{i=1}^n X_i$, and $N_t=\sup\{ n : \, S_n\leq t\}$. Now,
we denote by $\overline F (t)=P(X_1>t)$, the tail distribution function of the
$X_i$'s.

Let us now introduce a sequence $\eps_1, \eps_2,...$ of independent
symmetric Bernoulli variables i.e. satisfying $P(\eps_i=1)={1\over2}$ and
$P(\eps_i=-1)={1\over2}$.

The two following assertions hold:

(i) For every law of $X_1$, the process $Z_t:=\ds{\eps_{N_t}\over\o
F(t-S_{N_t})}$ is a martingale with respect to its natural filtration. 

(ii) The martingale $Z_t$ is normal if and only if $X_1$ follows the
distribution of $$-\ln U+{1\over 2 U^2}-{1\over2}$$ where $U$ denotes a uniform
variable on $[0,1]$.
\end{prop}

Proof: The proof of (i) is based on the fact that,
conditionally on the event $(t-S_{N_t}=x)$, the probability for $Z_t$  to jump
between
$t$ and $t+h$ is equal to $\ds-{\o F'(x)\over\o F(x)}h+o(h)$, and if a jump
occurs, its expectation is equal to
$\ds-{\eps_{N_t}\over \o F}$. On the other hand, when there is no jump, the
increment of the process, between $t$ and
$t+h$, is equal to
$\ds-\eps_{N_t}{\o F'(x)\over (\o F(x))^2}h+o(h)$. As a result, the expectations
of the increments in both situations are balanced. 

Let us now prove (ii). The martingale is normal if and only if the
conditional expectation of the square of the increment of $Z_t$ between $t$ and
$t+h$ is always equal to
$h+o(h)$. Again, conditionally on the event $(t-S_{N_t}=x)$, the probability for
$Z_t$ to jump between 
$t$ and $t+h$ is equal to $\ds-{\o F'(x)\over\o F(x)}h+o(h)$, and if a jump
occurs, the expectation of its square is equal to
$\ds{1\over2}((1+{1\over \o F})^2+(1-{1\over \o F})^2)$ which is $\ds 1+{1\over
\o F^2}$. Now, the absence of a jump contributes only in a $o(h)$ in the
expectation of the square of the increment of the process.

Therefore, the martingale is normal if and only if 
$$ -({1\over
\o F(x)}+{1\over
\o F^3(x)})\o F'(x)=1. $$
 This yields $ \ds -\ln \o F(x) +{1\over 2 \o F^2(x)} -{1\over 2}=x $, so that
the distribution function of $X_1$ is the inverse of $\ds x\mapsto -\ln(1-x)
+{1\over 2 (1-x)^2} -{1\over 2}$. This gives the result. \qed

{\it Remark:} We might have included in Proposition 1, analogous results
concerning the sometimes called ``second Az\'ema martingale" i.e. $\ds\sqrt{
\pi\over2}\sqrt{t-g_t}-l_t$, where
$l_t$ denotes the local time of zero at time $t$. Indeed, for every law of
$X_1$, 
$\tilde Z_t:=\ds{1\over\o F(t-S_{N_t})}-N_t$ is a martingale  which is normal if
and only if $X_1$ follows the distribution of $\ds{1\over 2 U^2}-{1\over2}$.  

But we shall not insist on that fact since this second Az\'ema martingale
is neither markovian nor enjoys the chaotic representation property, and does
not enter the class of Az\'ema martingales solving structure equations.

Let us denote by $Z^{(-1)}_t$ the normal martingale $Z_t$
of (ii).
We note that $Z_t^{(-1)}$ is a Markov process, and more
precisely,

\begin{prop}
The process $Z^{(-1)}_t$ is a Markov process with generator:
$$L^{(-1)}g(x)= {{1\over2}(g(-1)+g(1))-g(x)+xg'(x)\over 1+x^2}      $$
\end{prop}

Proof: From Proposition 1 (ii), we deduce that 
$$-{\ln \o F(x)}+{1\over 2\o
F^{2}(x)}-{1\over2} =x $$ 
which implies, by differentiation:
$$ -{\o F'(x)\over\o F(x)}={1\over1+\ds{1\over \o
F^{2}(x)}}$$ 
Hence, the rate of jump of the process at
time $t$, which is equal to  the value of the function $ \ds-{\o F'\over\o F}$
 at $(t-S_{N_t})$, is precisely ${1\over\ds 1+\ds (Z^{(-1)}_t)^2}$, and the
speed of the trajectory at a time $t$ between two jumps is equal to 
$\ds-\eps_{N_t}{\o
F'(t-S_{N_t})\over\o
F^2(t-S_{N_t})}$ which is precisely $\ds{Z^{(-1)}_t\over\ds 1+
(Z^{(-1)}_t)^2}$.
\qed

\begin{cor}
The process $Z^{(-1,n)}_t:={Z^{(-1)}_{nt}\over\sqrt n}$ is a Markov process with
generator:
$$ L^{(-1,n)}g(x)={{1\over2}(g(-{1\over \sqrt n})+g({1\over\sqrt
n}))-g(x)+xg'(x)\over {1\over n}+x^2}     
$$
\end{cor}
Proof: Consider a differentiable function $g$. Introduce the function
$h(x)=g({x\over\sqrt n})$. 

The image by $ L^{(-1,n)}$ of the function $g$ is given by:
$$ L^{(-1,n)}g=n \,.\,(L^{(-1)}h)(\sqrt n\,.\, x)$$ 
A direct computation gives the announced expression.\qed

This last corollary already gives an idea of the proximity between  $Z^n_t$ and
the first Az\'ema martingale. Indeed, as writes Emery in \cite{E1} about Az\'ema
martingales: ``formally (this means: informally !) it should be Markov, with
generator
$$Lg(x)=\left\{
\begin{array}{lc}
\ds{g((1+\beta)x)-g(x)-\beta xg'(x)\over \beta^2 x^2}&\hbox{if } x\neq0\\
{1\over2}g''(x)&\hbox{if } x=0"
\end{array}\right.$$
(in the case of the first Az\'ema martingale, the parameter $\beta$ is equal to
$-1$)

This is the aim of the next section which is to prove, among other things, the
convergence in distribution in the sense of the Skorohod topology of
${Z^{(-1)}_{nt}\over\sqrt n}$ towards the first Az\'ema martingale.

\section{ An invariance principle }

Let us begin with the definition of the  Az\'ema martingales. 

Following Emery, let us introduce a real parameter $\beta$. The structure
equation 
\begin{equation}
d[X,X]_t=dt+\beta X_{t^-}dX_t
\end{equation}
has a unique solution.

The existence was proven in a note of Meyer
\cite{M2}, and the uniqueness of the solution was proven by Emery in \cite{E1}
(for $\beta<0$) and \cite{E2} (for $\beta >0$).  The solution of this structure
equation is called the Az\'ema martingale with parameter $\beta$. 

When a jump occurs at time $t$, the value of the process
changes by some
prescribed factor.  
 More precisely, $X_t=(1+\beta)X_{t^-}$.

$\bullet$ The case
$\beta=0$  corresponds to Brownian motion: the continuous normal
martingale. 

$\bullet$ The case
$\beta=-1$ is solved by the martingale
$(sign\, B_t)\sqrt {2(t-g_t)} $ mentioned in the previous section. 

$\bullet$ The case
$\beta=-2$ is solved by Parthasarathy's martingale satisfying $|X_t|=\sqrt t$,
and who changes of sign according to a Poisson point process with intensity
$dt/4t$.

$\bullet$ In the case $\beta\leq -2$, the solution, at least starting
at a non-null point, is rather easy to define, since $|X_t|$  always goes away
from 0. In proving the uniqueness of the
solution, Emery  in \cite{E1} shows how to define the process starting from 0
by using  a self-similarity argument and a representation of the process by using
a time-changed Poisson process.

For $\beta>-2$, things go in a more complicated way since the process $X_t$
reaches 0 in finite time, and the above formal generator cannot be of much
help after that time. 

Let us mention finally that the problem of the chaotic representation property
still remains a challenging open question for parameters $\beta$ which do not
belong to
$[-2,0]$.

Now, let us turn to our main purpose concerning the statement of an invariance
principle for these processes. This question is quite natural, since
Az\'ema martingales are self-similar in the sense that, for all $\lambda>0$ the
processes $X_t$ and $X_{\lambda^2t}/\lambda$ are equal in distribution (one can
indeed easily check that the process $X_{\lambda^2t}/\lambda$ solves the right
structure equation, and conclude by uniqueness in law of the solution).

In order to define the process which will be, in our statement, at the origin of
such a limit theorem, we generalize the definition of $Z^{(-1)}_t$ introduced in
the previous section, and define for all real parameters $\beta$ the Markov
process $Z^{(\beta)}_t$ having generator:
$$ L^{(\beta)}= {{1\over2}(g((1+\beta)x-1)+g((1+\beta)x+1))-g(x)-\beta
xg'(x)\over 1+\beta^2x^2}   $$
Let us denote the process $Z^{(\beta,n)}_t:={Z^{(\beta)}_{nt}\over\sqrt n}$. It
has generator:
$$ L^{(\beta,n)}= {{1\over2}(g((1+\beta)x-{1\over \sqrt
n})+g((1+\beta)x+{1\over\sqrt n}))-g(x)-\beta xg'(x)\over {1\over
n}+\beta^2x^2}     
$$
We can now state our main result:

\begin{theo}
The sequence of processes $Z^{(\beta,n)}_t:={Z^{(\beta)}_{nt}\over\sqrt n}$
converges to the Az\'ema martingale with parameter $\beta$, in the sense of
the weak convergence for the Skorohod topology.
\end{theo}

Proof: The image by $ L^{(\beta,n)}$, of $x$ and $x^2$, being respectively $0$
and
$1$, we deduce that $Z^{(f,n)}_t$ are normal martingales. Therefore, by
Rebolledo's theorem (see \cite{R}, II.3.1), this sequence is tight for the weak
convergence in the Skorohod topology. All we have to prove is that all
limiting points of the sequence of processes $Z^{(\beta,n)}_t$ satisfy the
structure equation (1). The uniqueness in law of the solution of (1) allows to
conclude.

Let us define at a jump time $t$ of the process $Z^{(\beta,n)}_t$,
the symmetric Bernoulli variable $\eps^{(\beta,n)}_t$ defined by 
 $$Z^{(\beta,n)}_t-Z^{(\beta,n)}_{t^-}=\beta
Z^{(\beta,n)}_t+{\eps^{(\beta,n)}_t\over\sqrt n}$$
 Let us consider
$[Z^{(\beta,n)},Z^{(\beta,n)}]_t$. (We replace, in the sequel, the
superscript 
$(\beta,n)$ by $n$.)
$$\begin{array}{rl}[Z^{n},Z^{n}]_t&=\ds\sum_{s\leq t}(Z^n_s-Z^n_{s^-})^2= 
\sum_{s\leq t}(\beta Z^n_{s^-}+{\eps^n_s\over\sqrt n})(Z^n_s-Z^n_{s^-})\\
&=\ds\int_0^t \beta Z^n_{s^-}dZ_s^n-\int_0^t \beta Z^n_{s^-}1_{\Delta
Z_s^n\neq0}dZ_s^n +\sum_{s\leq t}1_{\Delta Z_s^n\neq0}{\eps^n_s\over\sqrt
n}(\beta Z^n_{s^-}+{\eps^n_s\over\sqrt n})\\
&=\ds\int_0^t \beta Z^n_{s^-}dZ_s^n+t-\int_0^t {ds\over
1+n(\beta Z^n_{s^-})^2}
+\sum_{s\leq t}1_{\Delta Z_s^n\neq0}{\eps^n_s\over\sqrt n}(\beta
Z^n_{s^-}+{\eps^n_s\over\sqrt n})
\end{array} $$

Our task, now, is to show that  the two last terms vanish when $n$ goes to
infinity, whereas the other terms converge to their analogous quantity for the
limiting process, so that every limiting process will solve the structure
equation satisfied by an Az\'ema martingale, and will therefore be identified
with the Az\'ema martingale with parameter $\beta$.

  Now we proceed, like in Meyer's proof of the existence of
solutions of structure equations, and introduce successively:

\_  a probability space $(\Omega, {\mathcal G}, Q)$, obtained from  the Skorohod
representation theorem, on which a subsequence
 of processes $X^{\sigma(n)}_s$, following the law of an extracted sequence of
$Z^n_s$  converges almost surely in $\lambda\otimes Q$. For ease of reading,
we keep writing $Z^n_s$ instead of $X^{\sigma(n)}_s$. In other words, we
 suppose that, for almost every time $t$, $Z^n_t$ converges almost surely.

\_ a sequence of stopping times $T_N^n=\inf\{t : |Z_t^n|>N\}$ where $N$ is large
enough, so that $Q(T_N^n\leq t)$ is smaller than any prescribed $\epsilon$,
uniformly in $n$ (Doob's inequality makes it possible to do so, since
the processes $Z_t^n$ have the same variance).

Let us  compare, in a first step, $[Z^{n},Z^{n}]_{T_N^n\wedge t}$ with the sum of
the square of its increments along some subdivision $(t_i)$ of the interval
$[0,t]$:
$$[Z^{n},Z^{n}]_{T_N^n\wedge t}-\sum_i(Z^n_{T_N^n\wedge t_{i+1}}-Z^n_{T_N^n\wedge
t_{i}})^2=2\sum_i\int_{t_i}^{t_{i+1}}(Z^n_{T_N^n\wedge
s_-}-Z^n_{T_N^n\wedge t_i})dZ^n_{T_N^n\wedge s}$$
The compensator of $Z_{T_N^n\wedge s}$ being equal to $T_N^n\wedge s$, we deduce,
$$E[([Z^{n},Z^{n}]_{T_N^n\wedge t}-\sum_i(Z^n_{T_N^n\wedge
t_{i+1}}-Z^n_{T_N^n\wedge t_{i}})^2)^2]\leq
4E[\sum_i\int_{t_i}^{t_{i+1}}(Z^n_{T_N^n\wedge s_-}-Z^n_{T_N^n\wedge t_i})^2ds]$$
which in turn is smaller that $4E[\sum_i\int_{t_i}^{t_{i+1}}(Z^n_{T_N^n\wedge
s_-}-Z^n_{T_N^n\wedge t_i})^2ds]\leq{4\over3}\sum_i(t_{i+1}-t_i)^3$.
(In this first step, stopping the martingales at $T_N^n$ just ensures the
integrability we need. In the further steps, this stopping argument will be
more widely used.) 

  As a conclusion, we deduce that the "discrete quadratic Riemann
sums" of $Z^{n}$ approximate  $[Z^n,Z^n]$ uniformly, in probability.  Now, since
the corresponding sums for $Z$ approximate $[Z,Z]$ in probability, it just
remains to notice that we can choose the time of our subdivisions among the
times $t$ for which $Z^{n}_t$ converges almost surely to $Z_t$ so that the
quadratic sums of $Z^n$ will converge almost surely and therefore in probability
to their corresponding analog for the process $Z$.

In a second step, we deal with the second term $\int_0^t Z^n_{s^-}dZ_s^n$. We
prove its convergence to the corresponding quantity for $Z$, exactly the same
way as in Meyer's proof.

Now, we have to prove the convergence to 0 in probability of both remaining
terms  $\ds\int_0^t {ds\over
1+n(\beta Z^n_{s^-})^2}
$ and $\ds\sum_{s\leq t}1_{\Delta Z_s^n\neq0}{\eps^n_s\over\sqrt n}(\beta
Z^n_{s^-}+{\eps^n_s\over\sqrt n}). $

Again, it suffices to prove it for these quantities when $t$ is replaced by 
$T_N^n\wedge t$. 

This result will be a consequence of the following key lemma, which gives some
uniform control (in $n$) on the time spent by the processes $Z_t^n$ near the
origin, and on the number of their jumps:

\begin{lem}

The two following convergence results hold: for all $N>0$, and all $t>0$,

(i) $\forall \epsilon>0,\, \exists \delta>0,\, \exists n_0\in\N,$ such that 
$\forall n\geq n_0$,
$$P(\int_0^{T_N^n\wedge t}1_{Z^n_s\in[-\delta,\delta]}ds>\epsilon)<\epsilon $$

(ii) $\ds{N^n_{T_N^n\wedge t}\over n}$ converges to 0 in probability.
\end{lem}

Proof: We start with (i). 
We will use the following notation for a constant which will appear repeatedly:
$C=\max(1,1+|\beta|)$.

Let us fix $\eps>0$, and, for some
$\delta\in ]0,{N\over 2C}[$, consider the successive stopping
times:

\_ $\tau^n_1=\inf\{t : |Z_t^n|>2\delta\}$

\_ $\tau^n_2=\inf\{t >\tau^n_1: |Z_t^n|<\delta \quad or \quad |Z_t^n|>N \}$

\_ If $\tau^n_2\neq T_N^n$, we define  $\tau^n_3:=\inf\{t  >\tau^n_2:
|Z_t^n|>2\delta\}$

and more generally for $i\geq2$,

\_ $\tau^n_{2i}:=\inf\{t >\tau^n_{2i-1}:
|Z_t^n|<\delta
\quad or
\quad |Z_t^n|>N \}$

\_ If $\tau^n_{2i}\neq T_N^n$, then $\tau^n_{2i+1}:=\inf\{t  >\tau^n_{2i}:
|Z_t^n|>2\delta\}$

Let us denote by $K$, the {\it random} integer $K:=\inf\{i\geq1: \tau^n_{2i}=
T_N^n\}$.

From these definitions, we get obviously:
$$\int_0^{T_N^n\wedge t}1_{Z^n_s\in[-\delta,\delta]}ds\leq
\tau^n_1+(\tau^n_3-\tau^n_2)+...+(\tau^n_{2K-1}-\tau^n_{2K-2})$$
Now, we can bound the expectation of the right hand term, by first noticing
that,  by the optional sampling theorem applied to the martingale $(Z^n_t)^2-t$, 
$$\begin{array}{rl}\forall i\geq 1,\quad
E[\tau^n_{2i+1}-\tau^n_{2i}\,\,|\,K>i]&=E[(Z^n_{\tau^n_{2i+1}})^2-
(Z^n_{\tau^n_{2i}})^2\,\,|\,K>i]\\
&\leq E[(Z^n_{\tau^n_{2i+1}})^2\,\,|\,K>i]\leq
4\max(1,(1+\beta)^2)\delta^2=4C^2\delta^2
\end{array}$$
 And secondly, we can estimate
$P(K=i+1\,\,|\,K>i)$, by exploiting that
$$E[Z^n_{\tau^n_{2i+2}}\,\,|\,K>i]=E[Z^n_{\tau^n_{2i+1}}\,\,|\,K>i] .$$

But, on one hand, $|Z^n_{\tau^n_{2i+1}}|  
\geq 2\delta$ and, on the other hand,   
$E[Z^n_{\tau^n_{2i+2}}\,\,|\,K>i]$ decomposes into
$E[1_{K=i+1}Z^n_{\tau^n_{2i+2}}]+ E[1_{K>i+1}Z^n_{\tau^n_{2i+2}}]$, with 
$$|E[1_{K=i+1}Z^n_{\tau^n_{2i+2}}]|\leq NC.P(K=i+1|K>i)   $$
and
$$|E[1_{K>i+1}Z^n_{\tau^n_{2i+2}}]|\leq \delta   $$ 
Therefore,
$$P(K=i+1\,\,|\,K>i)\geq {\delta\over NC}$$
which implies that  $ P(K\geq i) \leq (1-{\delta\over NC})^{i-1}$.

Now, writing
$$E[\tau^n_1+(\tau^n_3-\tau^n_2)+...+(\tau^n_{2K-1}-\tau^n_{2K-2})] 
=\sum_{i\geq1}E[(\tau^n_{2i-1}-\tau^n_{2i-2})1_{K\geq i}] $$
and, using the strong Markov property at time $\tau^n_{2i-1}$,  we get
$$E[\tau^n_1+(\tau^n_3-\tau^n_2)+...+(\tau^n_{2K-1}-\tau^n_{2K-2})]
\leq 4C^2\delta^2\sum_{i\geq1}(1-{\delta\over NC})^{i-1}=4NC^3\delta$$

We just have to choose $\delta$ such that $4NC^3\delta =\eps^2$ and
use Markov inequality, to get (i).

To prove (ii), we will strongly exploit (i). To make the reading easier, we will
fix $t$ equal to 1. Let us divide the interval
$[0,t]=[0,1]$ into $n$ intervals of the type $I_k^n=[{k\over n}, {k+1\over
n}]$, where $0\leq k\leq n-1$. 

Now, $\eps$ being fixed, we introduce $\delta$ like in (i), and define the
stopping times:
$$\sigma_k^n=\inf\{s\in I_k^n,\quad |Z_s^n|>\delta\} $$
with $\sigma_k^n=\infty$ if $\forall s\in I_k^n$, $|Z_s^n|\leq\delta$.

From (i), we get that, for $n$ large enough, with
probability bigger than
$1-\eps$, for more than
$(1-\eps) n$ integers
$k$, between $0$ and $n-1$, we have
$\sigma_k^n<\infty$. 

Let us introduce now the stopping time 
$\tau_k^n:=  \inf\{s>\sigma_k^n\quad |Z^n_s-Z^n_{s-}|\neq0\} $, with
$\tau_k^n=\infty$ whenever $\sigma_k^n=\infty$.

Now, for $n$  large enough, the quantity $P(|\tau_k^n-\sigma_k^n|<{2\over
n}\,|\,\sigma_k^n<\infty)$ can be bounded uniformly in $k$ by some sequence
depending on $n$, which  converges to 0. Indeed, if at some time $s$,
$|Z_s^n|>\delta$, and if moreover there is no jump between $s$ and $s+{2\over
n}$, then $|Z_s^n|$ remains bigger than $\delta/2$ for $n$ large enough. Now, on
the time interval
$[s,s+{2\over n}]$, the intensity of the jump remains smaller than  ${1\over
\beta^2(\delta/2)^2}$.

We deduce that, conditionally  on the event $(\sigma_k^n<\infty)$, 
$|\tau_k^n-\sigma_k^n|\wedge{2\over n}$ is stochastically bigger than 
$\e_{4\over\beta^2\delta^2}\wedge{2\over n}$, where
$\e_{4\over\beta^2\delta^2}$ is an exponential variable having parameter 
${4\over\beta^2\delta^2}$. 

So that, for $n$ large enough,
$$P(|\tau_k^n-\sigma_k^n|>{2\over
n}\,|\,\sigma_k^n<\infty)>\exp (- {8\over\beta^2\delta^2n})$$ 

Now,  $|\tau_k^n-\sigma_k^n|>{2\over
n}$ implies that there is no jump inside the time interval $I_{k+1}^n$.

So, we deduce that, for $n$ large enough, with
probability bigger than
$1-2\eps$, the process $Z_s^n$  will not jump on more than
$(1-2\eps) n$ intervals among the $I_k^n$'s. 

Finally, since the rate jump of the process is permanently less than $n$ and
the length of each $I_k^n$ is equal to ${1\over n}$,  for
$n$ large enough, we get that, with an arbitrarily high probability, there are
not more than $4\eps n$ jumps on the remaining $2\eps n$ intervals $I_k^n$ which
may have some jumps.
\qed  

{\it End of the proof of Theorem 1:}
The (i) part of Lemma 1 gives immediately the convergence in probability to 0,
of 
$\ds\int_0^{T_N^n\wedge t} {ds\over
1+n(\beta Z^n_{s^-})^2}$.

Now, $\ds\sum_{s\leq {T_N^n\wedge t}}1_{\Delta Z_s^n\neq0}{\eps^n_s\over\sqrt
n}(\beta Z^n_{s^-}+{\eps^n_s\over\sqrt n})= {N^n_{T_N^n\wedge t}\over
n}+\ds\beta\sum_{s\leq {T_N^n\wedge t}}1_{\Delta Z_s^n\neq0}{\eps^n_s\over\sqrt
n} Z^n_{s^-}$

Part (i) of Lemma 1 implies the convergence to 0 of the first term
of this sum, and we bound the
$L^2$-norm of the second term:
$$ E[(\sum_{s\leq {T_N^n\wedge t}}1_{\Delta Z_s^n\neq0}{\eps^n_s\over\sqrt
n} Z^n_{s^-})^2]={1\over n}E[\sum_{s\leq {T_N^n\wedge t}}1_{\Delta Z_s^n\neq0}
(Z^n_{s^-})^2]\leq N^2. E[{N^n_{T_N^n\wedge t}\over
n}] $$
Now, for all $\eps>0$,
$$E[{N^n_{T_N^n\wedge t}\over
n}]\leq \eps+E[{N^n_{T_N^n\wedge t}\over
n}.1_{({N^n_{T_N^n\wedge t}\over
n}>\eps)}]\leq\eps+E[({N^n_{T_N^n\wedge t}\over
n})^2]^{1\over2}P({N^n_{T_N^n\wedge t}\over
n}>\eps)^{1\over2} $$
Now, coupling the jump times of the process $Z_t^n$ with the points of a Poisson
process with intensity $n$, we are able to bound 
$E[({N^n_{T_N^n\wedge t}\over n})^2]$ uniformly in $n$, by $t^2$.

We conclude that $\ds\sum_{s\leq {T_N^n\wedge t}}1_{\Delta Z_s^n\neq0}{\eps^n_s\over\sqrt
n}(\beta Z^n_{s^-}+{\eps^n_s\over\sqrt n})$ converges to 0 in probability, and
obtain Theorem 1.
\qed

{\it Remark:} We might have replaced the randomness of the jumps,
having distribution ${1\over2}(\delta_{-1}+\delta_1)$ in the case of the process
$Z_t^{(\beta,1)}$, by any other one, having compact supported distribution with
mean 0 and variance 1. In other words, the result of Theorem 1 remains valid, if
we take, instead of $Z_t^{(\beta,1)}$, any process having generator
$$ {\left(\ds\int_\R g((1+\beta)x+y)\mu(dy)\right)-g(x)-\beta xg'(x)\over
1+\beta^2x^2}  $$
where $\mu$ is a compact supported measure, with mean 0 and variance 1. Indeed,
one can easily check that we still deal with normal martingales.

\section{An extension to structure equations}

Actually, Meyer is able to construct solutions of the more general structure
equation 
\begin{equation}
 d[X,X]_t=dt+f(X_{t^-})dX_t
\end{equation}
where $f$ is an arbitrary continuous function.

The strategy of Meyer consists in solving a discretized
structure equation, which leads to a cascade of second degree polynomial 
equations satisfied by the discretized increments, having always exactly two
solutions of opposite signs, which are chosen with respect to the (unique)
probability preserving the martingale property of the constructed discrete
process.

We want to show in this section how our approach can be used to exhibit a
solution of the structure equation (which may be different from Meyer's in
the non-uniqueness cases). Inspired by the processes of previous section, we
introduce a sequence of normal martingales
$Z^{(f,n)}_t$ whose limiting points in the sense of the weak convergence for the
Skorohod topology will satisfy the structure equation (2). 

However, although our approach has the advantage to give an invariance principle
in the self-similar case of previous section, nonetheless the construction of
structure equations for general continuous functions induces delicate
discussions on the time spent near the zeroes of the function $f$, and we did
not work them out for a general continuous function $f$.  So, our aim, in this
section, is more to examplify things, than to arrive to the most general
statement, which might be valid for general continuous functions, but we will
let this question open here.

We define  $Z^{(f,n)}_t$, as the  Markov
process starting at 0, and having generator:
$$ L^{(f,n)}g(x)={{1\over 2}\left(g(x+f(x)+{1\over \sqrt n})+g(x+f(x)-{1\over
\sqrt n})\right)-g(x)-f(x)g'(x)\over {1\over n}+f(x)^2}
$$

\begin{prop}
Suppose $f$ is a continuous function with isolated zeroes.  Suppose, in
addition, that, in the neighbourhood of each zero $x_j$ of $f$,
$f(x_j+h)=o(\sqrt h)$.

Then, the limiting points of the sequence of processes
$Z^{(f,n)}_t$, in the sense of the weak convergence for the Skorohod topology, 
satisfy the structure equation (2).
\end{prop}
Proof: 
Like in the previous section, we start by noticing that 
the image by
$ L^{(f,n)}$, of
$x$ and
$x^2$, is respectively
$0$ and
$1$.  Indeed,
$${1\over 2}\left((x+f(x)+{1\over \sqrt n})+(x+f(x)-{1\over
\sqrt n})\right)-x-f(x)=0$$
and
$${1\over 2}\left((x+f(x)+{1\over \sqrt n})^2+(x+f(x)-{1\over
\sqrt n})^2\right)-x^2-2xf(x)={1\over n}+f(x)^2   $$
We deduce that $Z^{(f,n)}_t$ are normal martingales. Therefore, by
Rebolledo's theorem (see \cite{R}, II.3.1), this sequence is tight for the weak
convergence in the Skorohod topology.

Now, we follow the proof of Theorem 1, and write: 
$$[Z^{n},Z^{n}]_t=\ds\int_0^t f(Z^n_{s^-})dZ_s^n+t-\int_0^t
{ds\over 1+nf( Z^n_{s^-})^2}
+\sum_{s\leq t}1_{\Delta Z_s^n\neq0}{\eps^n_s\over\sqrt n}(
f(Z^n_{s^-})+{\eps^n_s\over\sqrt n})$$

Applying again Skorohod representation theorem, we are able to prove the
convergence of the  terms $[Z^{n},Z^{n}]_t$ and $\int_0^t f(Z^n_{s^-})dZ_s^n$
to the analogous quantities for the limiting process.

Now, the key point is to check the analogous property to Part (i) of Lemma 1.

 In other words, we have to prove, that, for every $x_j\in [-N,N]$  (there are a
finite number of them),
$\forall \epsilon>0,\, \exists \delta_j>0,\, \exists n_0\in\N,$ such that 
$\forall n\geq n_0$,
$$P(\int_0^{T_N^n\wedge
t}1_{Z^n_s\in[x_j-\delta_j,x_j+\delta_j]}ds>\epsilon)<\epsilon
$$
We then introduce the same stopping times as in Lemma 1, and 
 the differences with the proof of Lemma 1 come out when we estimate
$E[\tau^n_{2i+1}-\tau^n_{2i}\,\,|\,K>i]$:
$$\begin{array}{rl}\forall i\geq 1,\,
E[\tau^n_{2i+1}-\tau^n_{2i}\,\,|\,K>i]&=E[(Z^n_{\tau^n_{2i+1}})^2-
(Z^n_{\tau^n_{2i}})^2\,\,|\,K>i]\\
&\leq E[(Z^n_{\tau^n_{2i+1}})^2\,\,|\,K>i]\leq
(2\delta_j+\sup_{[-2\delta_j,
2\delta_j]}|f(x)|)^2=o(\delta_j)
\end{array}$$
whereas $P(K=i+1|K>i)$ remains bigger than a constant times $\delta_j$.

As a result, $E[\tau^n_1+(\tau^n_3-\tau^n_2)+...+(\tau^n_{2K-1}-\tau^n_{2K-2})]$
is bounded by a constant times ${(2\delta_j+\sup_{[-2\delta_j,
2\delta_j]}|f(x)|)^2 \over \delta_j }=o(1)$. 

As in Lemma 1, we choose
$\delta_i$ such that this fraction is equal to $\eps^2$.

We finish the proof the same way as for Theorem 1.\qed

{\it Remark 1:} The framework of Proposition 3 includes the case of asymmetric
Az\'ema martingales corresponding to a function $f$ of the form
$ax1_{x\geq0}+bx1_{x\leq0}$. In this case,  Phan proved in
\cite{Phan} that the solution of the structure equation is unique. From this
result, we deduce that these processes are self-similar and
Proposition 3 provides, in this case, an invariance principle.

{\it Remark 2:} Our assumptions in Proposition 3 are obviously not optimal: for
instance, if $f$ is null, we are back to the case of a simple random walk, and
the convergence to the corresponding solution of the structure equation which is,
in this case, the Brownian motion is obviously valid. It would not be hard to
check that one can allow $f$ to be null on a locally finite number of intervals,
and satisfy the assumptions of Proposition 3 on the complement of the union of
these intervals. In this case, the term $\sum_{s\leq t}1_{\Delta
Z_s^n\neq0}{\eps^n_s\over\sqrt n}( f(Z^n_{s^-})+{\eps^n_s\over\sqrt n})$ will
not vanish in the limit and will ``alternatively" contribute with $\int_0^t
{f( Z^n_{s^-})^2\over {1\over n}+f( Z^n_{s^-})^2}ds$ to build a compensator
equal to $t$.

It is also very likely that the condition on  $f$ to be $o(\sqrt h)$
near its zeroes is not optimal (even if our attempts to relax it failed).  
\section{Construction of the process $Z_t^{(\beta,n)}$} 

We might just present the construction of $Z_t^{(\beta,1)}$ and rescale it, but
we prefered to keep the parameter $n$ in evidence, in order to point out
clearly where the randomization of the jumps appears. 

Let us introduce the functions $\o F$ and $H$ defined on $\R\times\R_+$, in
the following way: for any $(x_0,t)\in\R\times\R_+$, 

\_ $\o F(x_0,t)$ denotes the probability that the first jump of the process
$Z_t^{(\beta,n)}$, starting at $x_0$, takes place after time $t$.

\_ $H(x_0,t)$ is the value of $Z_t^{(\beta,n)}$, starting at $x_0$,
conditioned on the event that the first jump of the process takes place after
time $t$.

From the expression of the generator, we deduce the following system of
differential equations:

$$\left\{
\begin{array}{ccc}
\ds{1\over \o F}{\partial \o F\over\partial t}(x_0,t)&=&\ds-{1\over {1\over
n}+\beta^2 H^2(x_0,t)}\\ 
\ds{\partial H\over\partial t}(x_0,t)&=&\ds-{\beta H(x_0,t)\over {1\over
n}+\beta^2 H^2(x_0,t)}
\end{array}\right.$$

with the boundary conditions $\o F(x_0,0)=1$ and $H(x_0,0)=x_0$.

After noticing that $\ds{1\over\beta}{\partial \ln H\over\partial t}=
{\partial \ln \o F\over\partial t}$, we deduce that 
$$ H(x_0,t)=x_0(\o F(x_0,t))^\beta $$
From this fact, we deduce the following autonomous equation for $\o F$:
$$-{1\over n}{1\over \o F}{\partial \o F\over\partial t}-\beta^2x_0^2\o
F^{2\beta-1} {\partial \o F\over\partial t}=1
$$ 
which yields
$$-{1\over n}{\ln \o F}-{\beta x_0^2\over2}(\o
F^{2\beta}-1) =t $$ 

Therefore, the first jump of the process $Z_t^{(\beta,n)}$, starting
at
$x_0$, follows the law of $$-{1\over n}{\ln U}-{\beta x_0^2\over2}(
U^{2\beta}-1)$$
where $U$ is a uniform random variable on $[0,1]$. We notice here that the
special case $n=1$, $\beta=-1$ and $|x_0|=1$ allows us to recover the law of the
random variable we discussed in Proposition 1 (ii).

Denote $T_0=0$ and $T_i$ the $i$-th time of jump of the process
$Z_t^{(\beta,n)}$. We can describe the process $Z_t^{(\beta,n)}$ as follows:
there exists a sequence $(U_i)_{i\geq1}$ of independent uniform variables on
$[0,1]$, and a sequence $(\eps_i)_{i\geq1}$ of independent symmetric
Bernoulli variables, such that, for all $i\geq1$, (we omit here the superscript
$(\beta,n)$.)  
$$\left\{
\begin{array}{ccl}T_i&=&\ds T_{i-1}-{1\over n}\ln U_i-\beta {Z_{T_{i-1}}^2\over
2}(U_i^{2\beta}-1)\\
Z_{T_i^-}&=&\ds Z_{T_{i-1}} U_i^\beta\\
Z_{T_i}&=&\ds Z_{T_{i}^-}(1+\beta)+{1\over \sqrt
n}\eps_i
\end{array}\right.$$
On the time interval $[T_{i-1}, T_i[$, the trajectory $(t, Z_t)$ can be
parametrized by $(t(s), Z_{t(s)})_{0\leq s\leq 1-U_i}$ in the following way:
for all $s\in [0, 1-U_i[$,
$$\left\{
\begin{array}{ccl}t(s)&=&\ds T_{i-1}-{1\over n}\ln (1-s)-\beta
{Z_{T_{i-1}}^2\over 2}((1-s)^{2\beta}-1)\\
Z_{t(s)}&=&\ds Z_{T_{i-1}} (1-s)^\beta
\end{array}\right.$$

{\it Remark:} From this description of $Z_t^{(\beta,n)}$, we can easily deduce
a coupling between $Z_t^{(-1,n)}$ and $Z_t^{(\beta,n)}$, for all $n$ and all
parameters
$\beta$. But, unfortunately, this coupling cannot be transported at the level of
the corresponding Az\'ema martingales by merely making $n$ tend to infinity.  

This description allows us to draw some pictures of the process $Z_t^{(\beta,n)}$
for large values of $n$ and different values of the parameter $\beta$, having
therefore some good approximations of trajectories of Az\'ema martingales.
However, in this paper, we do not intend to quantify any error term.

We finish this paper with some simulations of Az\'ema martingales, which were
obtained by writing a code in MATLAB.
It basically draws
$Z_t^{(\beta,1)}$, asks for the number of arches one wants to appear in the
graph, and rescale the obtained function in order to get a function defined on
$[0,1]$. In the figures below, we asked for 1000 arches to appear.

Figure 1 represents Parthasarathy's martingale, whose graph is included inside
the parabola of equation $|y|=\sqrt x$. Figure 2 represents the first Az\'ema
martingale.
We can see on Figure 4 that, for very small values of the parameter $\beta$, the
Az\'ema martingale is a perturbation of the Brownian motion. It was
considered as a candidate in \cite{DP} to replace the Brownian motion
in the Black-Scholes model. We chose $\beta=-0.07$ which was the parameter
formerly estimated by Dritschel and Protter to modelise the variations of some
asset prices. 

We finally couple figures 3 and 5, which correspond respectively to the
parameters $-0.5$ and 1, and notice similarities in their graphs up to
some symmetry and horizontal affinity. Parthasarathy conjectured,
after noticing the duality between the two formal generators of these processes
that, up to a time change, the first  process was the time reversal of the
other (the conjecture was formulated, of course, for general positive values of
the parameter). This conjecture was actually proven in an unpublished part of
Phan's thesis \cite{Ph} by a nice and unusual argument involving polynomial test
functions.

\bigskip
\bigskip
\bigskip
\bigskip
\bigskip

\centerline{\epsfig{file=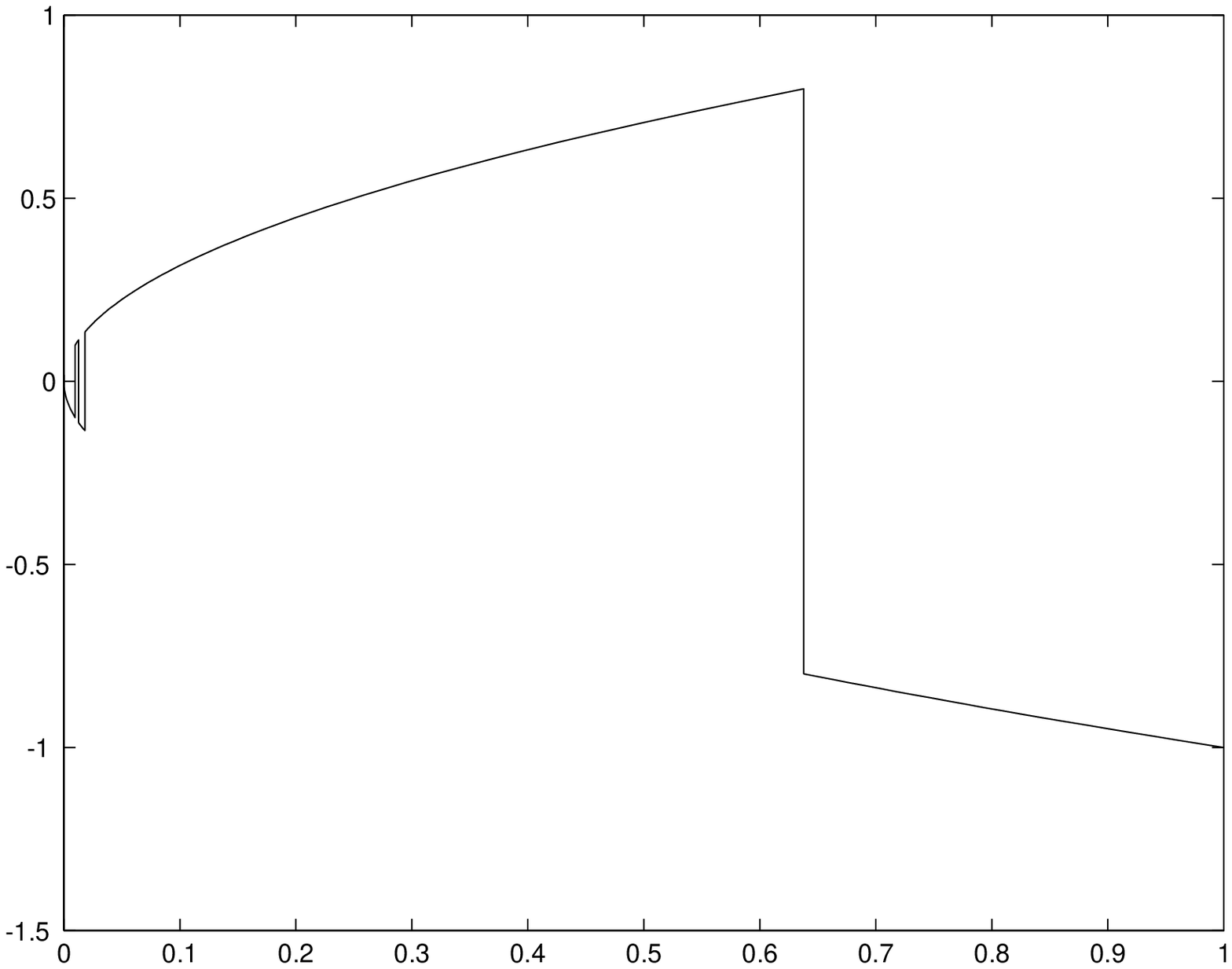, scale =0.36}}

\centerline{\small Fig. 1: $\beta =-2$}

\bigskip
\bigskip

\centerline{\epsfig{file=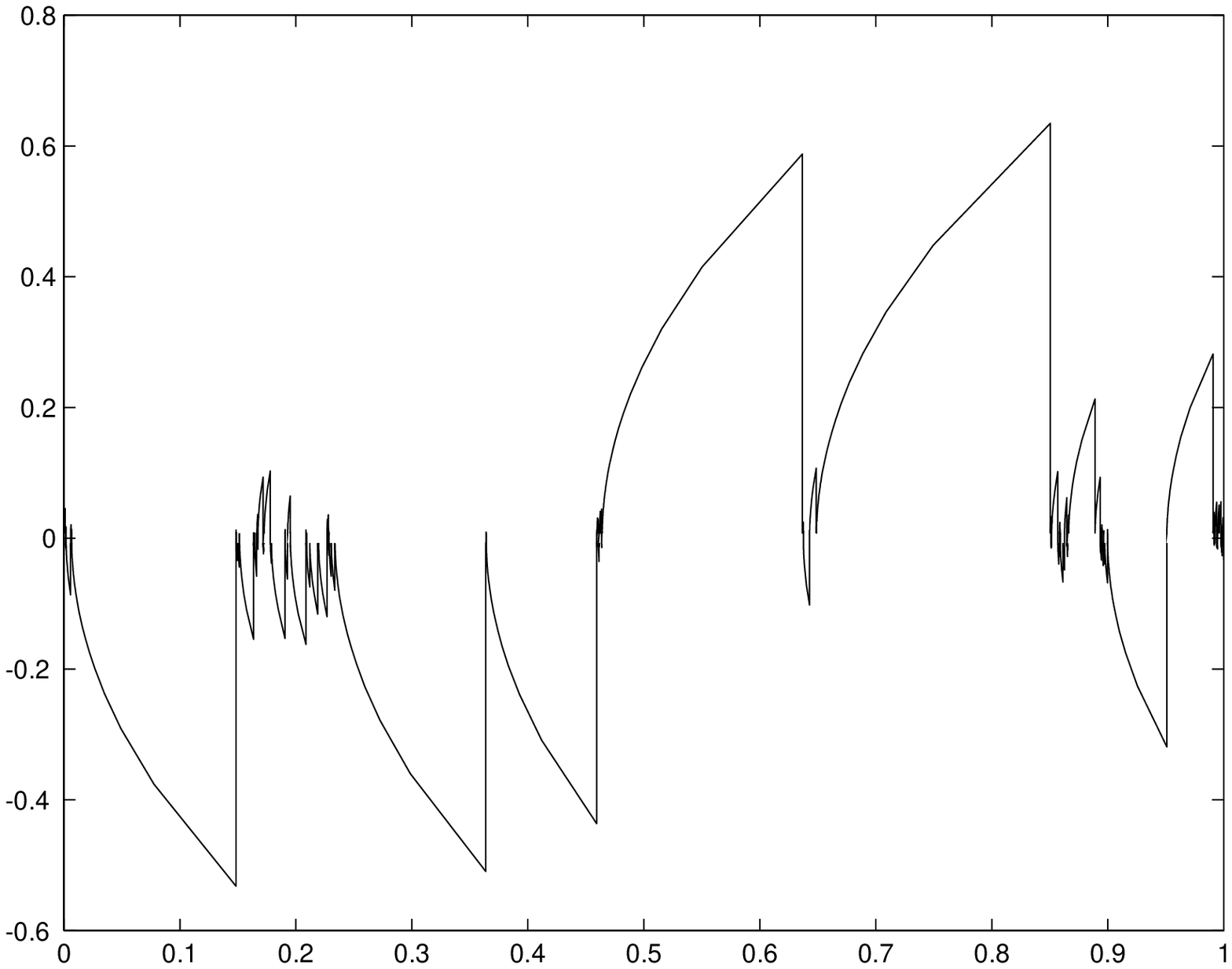, scale =0.36}}

\centerline{\small Fig. 2: $\beta =-1$}

\bigskip
\bigskip

\centerline{\epsfig{file=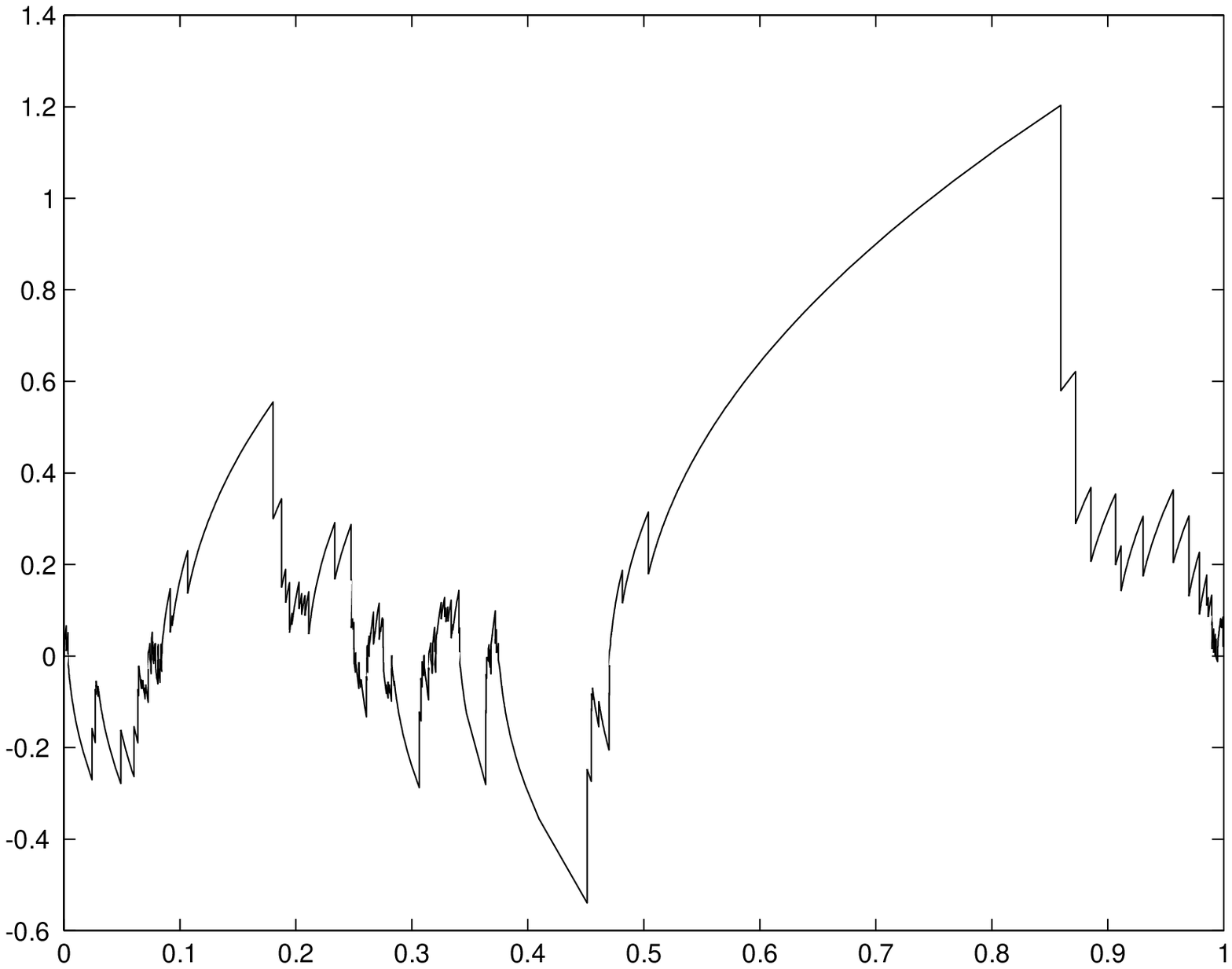, scale =0.36}}

\centerline{\small Fig. 3: $\beta =-0.5$}

\bigskip
\bigskip

\centerline{\epsfig{file=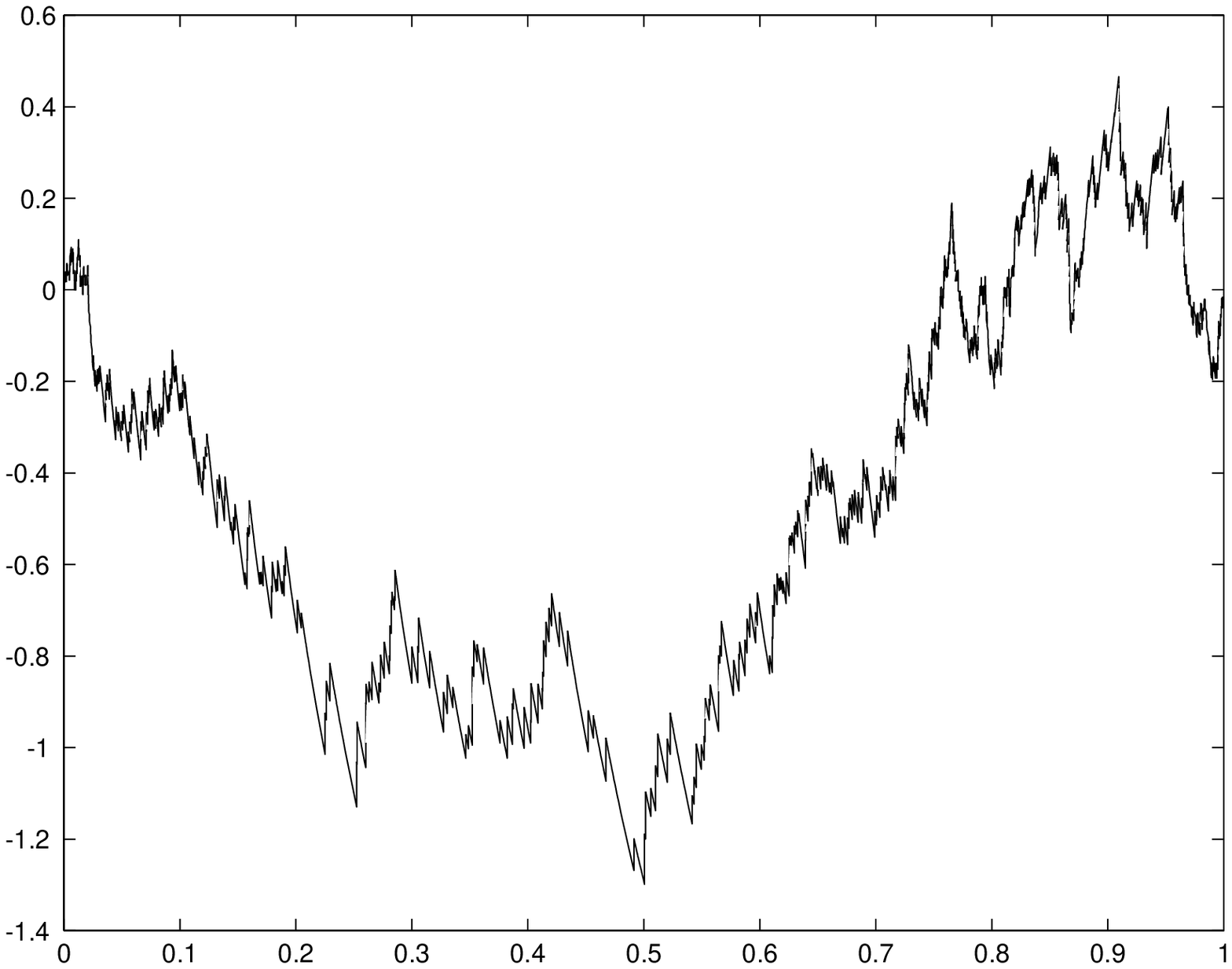, scale =0.36}}

\centerline{\small Fig. 4: $\beta =-0.07$}

\bigskip
\bigskip

\centerline{\epsfig{file=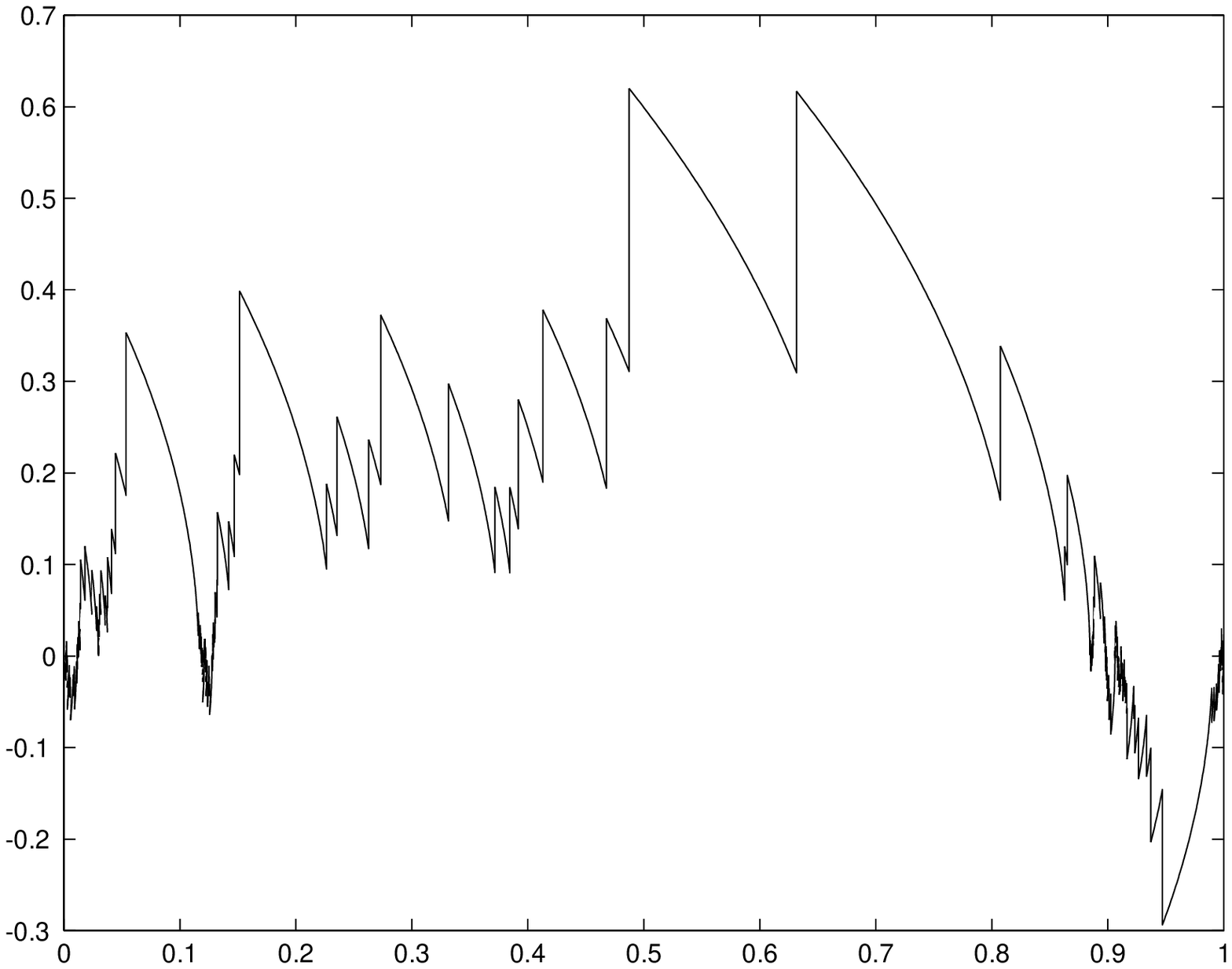, scale =0.36}}

\centerline{\small Fig. 5: $\beta =1$}

\end{document}